\documentclass{amsart}
\usepackage{graphicx}

\theoremstyle{definition}

\theoremstyle{remark}

\numberwithin{equation}{section} 

\begin{document}

\title[]{Complex centers of polynomial differential equations}%
\author{M.A.M. Alwash}%
\address{Department of Mathematics, West Los Angeles College\\ 9000 0verland Avenue, Los Angeles, CA 90230-3519, USA}%
\email{alwashm@wlac.edu}%

\thanks{I would like to thank Dr. Colin Christopher for a very fruitful email discussion.
I am very grateful to the Simmons Hall residential scholar program at MIT for their hospitality}%

\begin{abstract}
We present some results on the existence and nonexistence of centers for polynomial first order ordinary differential
equations with complex coefficients. In particular, we show that binomial differential equations without linear terms
do not have complex centers. Classes of polynomial differential equations, with more than two terms, are presented that
do not have complex centers. We also study the relation between complex centers and the Pugh problem. An algorithm is
described to solve the Pugh problem for equations without complex centers. The method of proof involves phase plane
analysis of the polar equations and a local study of periodic solutions.
\end{abstract}
\maketitle 
\section{Introduction}
Consider the differential equation
\begin{equation}
\dot{z}:=\frac{dz}{dt}=A_{N}(t) z^{N}+A_{N-1}(t)z^{N-1}+\ldots+A_{1}(t)z
\end{equation}
where $z$ is complex and $A_{i}(t)$ are continuous functions in $t$. Let $z(t,c)$ be the solution of (1.1) such that
$z(0,c)=c$. For a fixed real number $\omega$, we say that $z(t,c)$ is \emph{periodic} when $z(0,c)=z(\omega,c)$. If the
functions $A_{i}(t)$ are periodic with period $\omega$, then a periodic solution of equation (1.1) is a periodic
function with period $\omega$. The \emph{multiplicity} of a periodic solution $\varphi(t)$ of (1.1) is the multiplicity
of $\varphi(0)$ as a zero of the displacement function $q\longmapsto z(\omega,c)-c$. Note that $q$ is defined and
analytic in an open set containing the origin. The solution $z=0$ is called a \emph{center} for the differential
equation if all solutions starting in a neighborhood of $0$ are $\omega-$periodic. If the coefficients are real
functions and $A_{N}(t)$ does not change sign, then $z=0$ is not a center (see \cite{A1}). There are real centers when
$A_{N}(t)$ changes sign (see \cite{A4}, \cite{A5}, and \cite{A15}); the problem is related to the classical
center-focus problem of polynomial two-dimensional systems. The case $N=3$, with $A_{3}(t) \equiv 1$, was considered
recently in \cite{A6}; where it was shown that the equation could have a center at the origin if the coefficients are
complex valued. The work in \cite{A6} is mainly motivated by a problem stated in \cite{A8}. It was shown in \cite{A8}
that if the coefficients are real and bounded by a constant $C$, then there are at most
$8\exp\{(3C+2)\exp[\frac{3}{2}(2C+3)^{n}]\}$ periodic solutions. The problem in \cite{A8} is to prove an analogue of
this result for complex equations, that is, for complex $z$ and complex coefficients $A_{i}$. This paper is partly
motivated by the results of \cite{A6}.

Another motivation for the work in this paper is the local study of the Pugh problem about equations with real
coefficients.
 We recall that Pugh problem (see \cite{A13}) is to
find an upper bound for the number of periodic solutions in terms of $N$. It was shown in \cite{A9} that there are no
upper bounds for the number of periodic solutions for $N \geq 4$. Upper bounds can be found for some particular
classes; (see \cite{A1} and its references). Hence, one should seek upper bounds for the number of periodic solutions
in terms of $N$ and the degrees of the real polynomial functions $A_{i}(t)$. A local version of Pugh problem is to find
an upper bound for the multiplicity of a periodic solution in terms of $N$ and the degrees of the polynomials
$A_{i}(t)$. This problem was considered in \cite{A12}, with $N=4$. It was conjectured in \cite{A12}, that if
$A_{4}\equiv 1$, $A_{1}\equiv0$, $A_{2}$ and $A_{3}$ are of degree $k$ then the multiplicity of the origin is at most
$k+3$. It was shown in \cite{A3}, that when $k=2$ then the multiplicity of the origin is at most $8$ and there is a
unique equation with this maximum multiplicity. It was shown in \cite{A2}, that the multiplicity of the origin is at
most $10$ when the degrees of $A_{2}(t)$ and $A_{3}(t)$ are $2$ and $3$, respectively. Having determined the maximum
multiplicity, the next step is to construct equations with this number of periodic solutions. This is done by making a
sequence of perturbations in $A_{2}(t)$ and $A_{3}(t)$, each of which reduces the multiplicity of the origin by one; a
periodic solution thus bifurcates out of the origin. This bifurcation task was considered in \cite{A2} and \cite{A3}.
The local problem can be considered with the use of Groebner bases method. This task is considered in the last section;
the case $N=4$ and $A_{1}\equiv 0$ is studied. The problem reduces to study the solvability of system of polynomial
equations in many variables. Since the solvability in the theory of Groebner bases is over the field of complex
numbers, it becomes necessary to consider equations with complex coefficients. That is to consider complex centers.

All the known equations with complex centers have linear terms. On the other hand, computations for equations without
linear part demonstrate that $z=0$ is not a center for polynomial coefficients. These remarks lead us to conjecture
that polynomial differential equations without linear terms do not have centers at the origin, at least when the
coefficients are polynomial functions in $t$.
\newtheorem*{r2}{Conjecture}
\begin{r2}
Assume that $A_{j}(t)$, for $j=2,3,\cdots,N-1$, are polynomial functions. The solution $z=0$ is not a center for the
differential equation
\[ \dot{z}=z^{N}+A_{N-1}(t) z^{N-1}+\cdots+A_{3}(t) z^{3} + A_{2}(t)
z^{2} .\]
\end{r2}

First, we present classes of equations with complex coefficients that do not have centers at the origin. Our main
results in this direction are the following:
\newtheorem*{T1}{Theorem A}
\begin{T1}
Assume that the differential equation (1.1) has an invariant line $\theta = \alpha$, and either $A_{N}(t) > \gamma >
0$, or $A_{N}(t) < \gamma < 0$. If $\gamma \cos((N-1) \alpha) > 0$ then the equation can have only a finite number of
periodic solutions on the invariant line. In particular, $z=0$ is not a center.
\end{T1}
\newtheorem*{C1}{Corollary A}
\begin{C1}
Assume that $A_{N}(t) > \gamma >0$.
\begin{enumerate}
\item
 Let $A_{j}(t)=a_{j}(t)+ib_{j}(t)$, where
$a_{j}$ and $b_{j}$ are real continuous functions. If there exists an integer $k$ such that
\[ a_{j}\sin(j \frac{2k\pi}{N-1})=b_{j} \cos(j \frac{2k\pi}{N-1})\]
for $j=1,2,\cdots,N-1$, then $z=0$ is not a center for (1.1).
\item
If the coefficients $A_{j}$, $j=1,2,\cdots,N-1$, are real functions then $z=0$ is not a center for (1.1).
\item
 If $N-1$ is a multiple of $4$, $A_{j}$ are pure imaginary
functions for all odd $j$, and $A_{j}$ are real for all even $j$, then $z=0$ is not a center for equation (1.1).
\end{enumerate}
\end{C1}

\newtheorem*{T2}{Theorem B}
\begin{T2}
Consider continuous functions $A(t), B(t)$ and $C(t)$ and let
\[ \lambda_{1}=\int_{0}^{\omega} C(t) dt,\]
\[ \lambda_{2}=\int_{0}^{\omega} B(t) dt, \]
and
\[ \lambda_{3}=\int_{0}^{\omega}[A(t)+(M-L)B(t) (\int_{0}^{t} C(s) ds) ]dt . \]
The solution $z=0$ of the differential equation
\begin{equation}
\dot{z}=A(t) z^{N}+B(t) z^{M} +C(t) z^{L},
\end{equation}
with $1<L<M<N=L+M-1$, has multiplicity $L$ if and only if $\lambda_{1} \neq 0$. The multiplicity is $M$ if and only if
$\lambda_{1}=0$ and $\lambda_{2} \neq 0$; and it is $N$ if and only if $\lambda_{1}=\lambda_{2}=0$ and $\lambda_{3}
\neq 0$. If the origin is a center then $\lambda_{1}=\lambda_{2}=\lambda_{3} =0$.
\end{T2}
\newtheorem*{C2}{Corollary B}
\begin{C2}
\begin{enumerate}
\item
For any continuous functions $A(t)$ and $B(t)$ with $\int_{0}^{\omega} A(t) dt \neq 0$, the origin is not a center for
the equation
\[ \dot{z}=A(t) z^{N}+B(t) z^{M} \]
when $1<M<N$.
\item
Assume that $A(t), B(t), C(t)$, and $D(t)$ are continuous functions, and consider the differential equation
\[ \dot{z}=A(t) z^{N} + B(t) z^{M} + C(t) z^{L} + D(t) z ,\]
with $1<L<M<N=L+M-1$.
 Let
\[ A_{1}(t)=e^{(1-N) \int_{0}^{t} D(s) ds} A(t)\]
\[ B_{1}(t)=e^{(1-M) \int_{0}^{t} D(s) ds} B(t)\]
\[ C_{1}(t)=e^{(1-L) \int_{0}^{t} C(s) ds} C(t). \]
If any of the following conditions is not satisfied then the origin is not a center,
\[ \int_{0}^{\omega}C_{1}(t) dt=0,\]
\[ \int_{0}^{\omega}B_{1} (t) dt =0, \]
and
\[ \int_{0}^{\omega}[A_{1}(t) +(M-L) B_{1}(t)\int_{0}^{t} C_{1}(s)
ds]dt = 0.\]
\end{enumerate}
\end{C2}
\newtheorem*{r5}{Remark}
\begin{r5}
It follows from Theorem B, that:
\begin{enumerate}
\item
If $B(t)$ and $C(t)$ are polynomial functions with small coefficients, then $z=0$ is not a center for the differential
equation
\[ \dot{z} = z^{N}+ B(t) z^{M} + C(t) z^{L} \]
with $1<L<M<N=L+M-1$.
\item
A polynomial differential equation with only two terms has a center at $z=0$ only when one of the terms is linear.
\end{enumerate}
\end{r5}

 Now, we describe some equations that have centers at the origin.
\newtheorem*{T3}{Theorem C}
\begin{T3}
The solution $z=0$ of the differential equation
\begin{equation}
\dot{z}=z^{N}+A(t)z
\end{equation}
is of multiplicity $N$ if and only if
\[e^{\int_{0}^{\omega} A(t) dt}=1,\,\,\int_{0}^{\omega}
e^{(N-1)\int_{0}^{t} A(s) ds} dt \neq 0.\] The origin is a center if and only if
\[e^{\int_{0}^{\omega} A(t) dt}=1,\,\, \int_{0}^{\omega} e^{(N-1)
\int_{0}^{t} A(s) ds} dt =0. \]
\end{T3}

Equations with centers at the origin are given in the following result; the second part is a version of the result in
\cite{A6}. As it was first noticed in \cite{A6}, the center variety in each of the equations considered here contains
infinitely many connected components. Center variety of each of the known real centers contains a finite number of
connected components.

\newtheorem*{C4}{Corollary C}
\begin{C4}
\begin{enumerate}
\item
 With $\omega=2 \pi$, we consider the differential equation
\[ \dot{z}=z^{N}+(C+\frac{2}{N-1} \tan(t+c_{0}))z, \]
where $c_{0}$ is a complex non-real number. The multiplicity of $z=0$ is one if and only if $C=pi$, where $p$ is an
integer. The multiplicity is $N$ if and only if
\[ \int_{0}^{2 \pi} e^{(N-1)pit} \sec^{2}(t+c_{0}) dt \neq 0.\]
If this integral vanishes then $z=0$ is a center. In particular, there is a center when $(N-1)p$ is an odd number.
\item
 Suppose that $C$ is a number and $x$ is a real number.
With $\omega=2\pi$, the solution $z=0$ is a center for the equation
\[ \dot{z}=z^{N}+\frac{1}{N-1}
(C-xi \cos t)z \] if and only if $C=pi$, where $p$ is an integer and $x$ is a zero of the Bessel's function $J_{p}$.
\end{enumerate}
\end{C4}
Now, we give another class of equations that have a center. The coefficients of the equation satisfy a composition type
condition similar to those in \cite{A4}, \cite{A5}, and \cite{A15}.
\newtheorem*{t4}{Theorem D}
\begin{t4}
 Let $p(t)$ be a
differentiable $\omega-$periodic function and $q(t)$ is a continuous function with $e^{\int_{0}^{\omega} q(t)dt} = 1$.
The solution $z=0$ is a center for the differential equation
\[ \dot{z}=q(t)\, z + \sum_{k=2}^{N} p'(t)\, f_{k}(p(t))\, e^{(1-k) \int_{0}^{t} q(s) ds}\, z^{k} ,\]
where $f_{k}$ are any continuous functions.
\end{t4}
Finally, we use the method of groebner bases to study multiplicity of periodic solutions when the coefficients are
polynomial functions in $t$, or in $\cos t$ and $\sin t$.
\newtheorem*{T4}{Theorem E}
\begin{T4}
Consider the following equation
\[\dot{z}= z^{4} + A(t) z^{3} + B(t) z^{2}, \]
where $A(t)$ and $B(t)$ are polynomial functions.\\ (I) In each of the following classes of coefficients, $z=0$ is not
a center and the maximum possible multiplicity is the same whether the coefficients are complex or restricted to be
real.
\begin{enumerate}
  \item $A$ and $B$ are polynomial functions in $t$ of degrees $3$ and $2$, respectively..
  \item $A$ and $B$ are polynomial functions in $t$ of degrees $5$ and $1$, respectively.
  \item $A$ and $B$ are polynomial functions in $t$ of degrees $1$ and $3$, respectively.
\end{enumerate}
(II) In each of the following classes of coefficients, $z=0$ is not a center and the maximum possible multiplicity for
real coefficients is less than the maximum multiplicity for complex coefficients. Moreover, the origin is not a
center.
\begin{enumerate}
\item $A$ and $B$ are homogeneous polynomial functions in $\cos t$ and $\sin t$ of degrees $1$.
\item $A$ and $B$ are homogeneous polynomial functions in $\cos t$ and $\sin t$ of degrees $2$.
\item $A$ and $B$ are homogeneous polynomial functions in $\cos t$ and $\sin t$ of degrees $3$ and $1$, respectivly.
\end{enumerate}
\end{T4}
 In the next section, we prove Theorem A. Section 3 contains the proof of Theorem B. The proofs of Theorem C and
 Theorem D
are given in Section 4. In Section 5, we consider the case $N=4$. Several classes of equations are considered and
provide evidences to support our conjecture.
\section{Proof of Theorem A}
\begin{proof}
Let $A_{m}(t)=a_{m}(t)+ib_{m}(t)$ and $z=r e^{i\theta}$. With this notation, equation (1.1) becomes
\[ (\dot{r}+ i \dot{\theta} r)e^{i\theta}=
r^{N}e^{Ni\theta}P_{N}+r^{N-1}e^{(N-1)i\theta}(a_{N-1}+ib_{N-1})+ +\cdots+r e^{i\theta}(a_{1}+ib_{1}).\] Multiplying
both sides by $e^{i\theta}$ and separating the real and imaginary parts, give:
\[ \dot{r}=r^{N}P_{N} \cos((N-1)\theta)+
\sum_{j=1}^{N-2}r^{N-j}[a_{N-j}\cos((N-j-1)\theta)-b_{N-j} \sin((N-j-1)\theta)] \]
\[ \dot{\theta}=r^{N-1} P_{N} \sin((N-1)\theta)+
\sum_{j=1}^{N-2}r^{N-j-1}[a_{N-j} \sin((N-j-1)\theta)+b_{N-j} \cos((N-j-1)\theta)]. \]

Now, let us choose a real number $a$ with
\[ a >(\frac{2}{(N-1)\gamma \omega \cos((N-1)\alpha)})^{\frac{1}{N-1}}\]
 and such that the inequality
\[ \dot{r}> \frac{\gamma \cos((N-1) \alpha)}{2}r^{N}\]
is satisfied when $r \geq a$ and $\theta=\alpha$. Consider a solution, $re^{i\theta}$ with initial condition
$r_{0}e^{i\alpha}$. If $r_{0} \geq a$, then the inequality implies that
\[ r(t) \geq \frac{1}{(r_{0}^{1-N}-0.5(N-1)\gamma \cos((N-1)\alpha) t)^{\frac{1}{N-1}}}.\]
Therefore, if $r_{0} \geq a$, then the solution becomes unbounded at a point $t_{0}$ with
 \[t_{0} \leq \frac{2 r_{0}^{1-N}}{(N-1)\gamma
\cos((N-1)\alpha)} \leq \omega .\]
 Hence, the solution is undefined
on the interval $[0,\omega]$. Therefore, solutions start outside the disk $r \leq a$ are not periodic, and any solution
leaves the disk $r \leq a$ stays outside the disk as time increases.

If there is an infinite sequence of periodic solutions $z(t,r_{n}e^{i \alpha})$. The numbers $r_{n}$ are inside the
disk $r \leq a$. Hence, let $r_{n}\rightarrow r_{0}$. If $z(t,r_{0}e^{i\alpha})$ is a periodic solution, then all
solutions starting in a neighborhood of $r_{0}e^{i\alpha}$ are periodic. But $q$ is an analytic function. It follows
that $q \equiv 0$. Therefore, $z=0$ is a center for the equation (1.2). We define a real number $R$ by
\[R=sup\{r:z(t,re^{i\alpha}) \,\,is\,\, periodic\}. \]
It is clear that $R<\infty$. If the solution $z(t,Re^{i\alpha})$ is defined, then it is a periodic solution. Since $q$
is analytic, it follows from the property of continuous dependence of a solution on its initial value, that
$z(t,re^{i\alpha})$ is periodic for $r$ in a neighborhood of $R$. This is contrary to the definition of $R$. Therefore,
$z(t,Re^{i\alpha})$ is undefined on the interval $[0,\omega]$ and becomes unbounded at $t_{0}\varepsilon [0,\omega]$
and leaves the disk $r \leq a$. But $z(t,Re^{i\alpha})$ leaves $r \leq a$, whence $z(t,re^{i\alpha})$ leaves $r \leq a$
for $r < R$ and close to $R$. These solutions will not return to the disk $r \leq a$ because $\dot{r}(t)>0$, and hence
they are not periodic solutions. This is a contradiction to the assumption that $z=0$ is a center.

If $z(t,r_{0}e^{i\alpha})$ is undefined, then it leaves the disk $r \leq a$. Therefore, $z(t,r_{n}e^{i\alpha})$ leaves
$r \leq a$ for large enough $n$. This is again a contradiction to the periodicity of $z(t,r_{n}e^{i\alpha})$.
\end{proof}

To prove the first part of Corollary A, let $\alpha =\frac{2k \pi}{N-1}$. The condition in the statement implies that
\[a_{j} \sin(j\alpha)=b_{j} \cos(j\alpha).\]
\[\sin((N-j-1)\alpha)=-\sin(\alpha),\,\,
\sin((N-1)\alpha)=0.\]
 Hence,
\[\dot{\theta}(\alpha)=0.\]
The line $\theta = \alpha$ is an invariant line. Any solution has a point on this line stays on the line as long as it
is defined.

 On the other hand,
\[ \cos((N-1)\alpha)=1.\]
Hence, $A_{N}(t) \cos((N-1)\alpha) = A_{N}(t) > \gamma$. The conditions of Theorem A are satisfied and the result
follows.

 The second part of
Corollary A, follows directly from the first part with $k=0$. It is known that $z=0$ is not a center when the
coefficients are real functions (see \cite{A1}). To prove the third part, we take $k=\frac{N-1}{4}$. This implies that
$\alpha =\frac{2k \pi}{N-1}=\frac{\pi}{2}$. Thus $\sin(m\alpha)=0$ if $m$ is even, and $\cos(m\alpha)=0$ if $m$ is
odd.

\newtheorem*{r1}{Remark}
\begin{r1}
Equations (1.1), with real coefficients, have been studied in an interesting paper of Lloyd \cite{A10} using the
methods of complex analysis and topological dynamics. However, the results of \cite{A10} do not hold when the
coefficients are allowed to be complex. For example, Abel differential equation, $N=3$, has $3$ periodic solutions when
its coefficients are real. On the other hand, it may have an infinite number of periodic solutions when the
coefficients are complex (see, \cite{A6}). However, the method of \cite{A10} could be modified to obtain global results
for equations with complex coefficients. This exploitation is more subtle than the real case, and we defer this to
another paper.
\end{r1}
\section{Proof of Theorem B}
\begin{proof}
For $0 \leq t \leq \omega$ and $c$ in a neighborhood of $0$, we write
\[z(t,c)=\sum_{n=1}^{\infty} d_{n}(t)c^n, \]
where $d_{1}(0)=1$ and $d_{n}(0)=0$ if $n>1$. Thus
\[q(c)=(d_{1}(\omega)-1)c+\sum_{n=2}^{\infty}d_{n}(\omega)c^{n}.\]
The multiplicity is $K$ if and only if
\[d_{1}(\omega)=1,d_{2}(\omega)=d_{3}(\omega)=\cdots=d_{K-1}(\omega)=0,d_{K}(\omega)\neq 0.\]
The origin is a center when $d_{1}(\omega)=1$ and $d_{n}(\omega)=0$ for all $n>1$. The functions $d_{n}(t)$ are
determined by substituting the sum into the equation (1.2) and comparing coefficients of powers of $c$. The following
recursive sequence of differential equations is obtained with the initial equations.
\[\dot{d}_{n}=A S_{N} +B S_{M} + C S_{L},\]
where
\[S_{K}=\sum_{j_{1}+j_{2}+\cdots+j_{K}=n}^{} d_{j_{1}}d_{j_{2}}\cdots d_{j_{K}},\]
with $K \epsilon \{N, M, L\}$. To solve these equations, we integrate repeatedly. It is clear that $\dot{d}_{1}=0$ and
hence $d_{1}\equiv 1$. The next non-zero equation is
\[\dot{d}_{L}=C (d_{1})^{L}.\]It gives
\[ d_{L}(t)=\int_{0}^{t} C(s)ds .\]
The formula for $d_{M}$ has two possibilities. If $M+1\neq 2L$, then
\[\dot{d}_{M}=B (d_{1})^{M}.\] This implies that
\[ d_{M}(t)=\int_{0}^{t} B(s) ds.\]
In the case $M+1=2L$, the equation becomes
\[\dot{d}_{M}=B(d_{1})^{M}+L C d_{L} (d_{1})^{L-1}).\] We integrate
this equation to obtain
\[d_{M}(t)=\int_{0}^{t} B(s) ds + \frac{L}{2} (\int_{0}^{t} C(s)
ds)^2.\] Consequently, the formula for $d_{N}$ has two possibilities. The equation is
\[ \dot{d}_{N}=A (d_{1})^{N}+B M d_{L}(d_{1})^{M-1}+C L
d_{M}(d_{1})^{L-1}.\] We integrate this equation to obtain $d_{N}(t)$. If $M+1 \neq 2L$, then
\[d_{N}(t)=\int_{0}^{t}A(s) ds+(M-L) \int_{0}^{t} (B(s) \int_{0}^{s} C(u) du)ds + L
\int_{0}^{t}B(s) ds \int_{0}^{t} C(s) ds.\] In the case $M+1=2L$, the formula becomes
\begin{gather*}
 d_{N}(t)=\int_{0}^{t} A(s) ds+(M-L) \int_{0}^{t} (B(s) \int_{0}^{s} C(u) du)ds\\
+L \int_{0}^{t} B(s) ds \int_{0}^{t} C(s) ds + \frac{L^{2}}{6} (\int_{0}^{t} C(s) ds)^{3} .
\end{gather*}
The multiplicity is $L$ if $\int_{0}^{\omega} C(t) dt \neq 0$. The multiplicity is $M$ if $\int_{0}^{\omega} C(t) dt
=0$ but $\int_{0}^{\omega} B(t) dt \neq 0$. Finally, if the multiplicity is greater than $M$, then
\[ d_{N}(\omega)=\int_{0}^{\omega}A(t) dt + (M-L) \int_{0}^{\omega} (B(t) \int_{0}^{t}
C(s) ds)dt.\]
 The assumption in the statement of Theorem B, implies
that $d_{N}(\omega) \neq 0$. Therefore, the multiplicity is at most $N$.
\end{proof}
The first part of Corollary B follows directly from the above result. If $B(t)\equiv 0$, then $d_{N}(\omega)> \gamma
\omega \neq 0$. Hence, the multiplicity is at most $N$.

Now we prove the second part of Corollary B. The first necessary condition for a center is $e^{\int_{0}^{\omega} D(t)
dt} = 1$. We make the transformation $w=e^{\int_{0}^{t} D(s) ds} z$, and obtain
\[ \dot{w}= A_{1}(t) w^{N} + B_{1} (t) w^{M} + C_{1} (t) w^{L}, \]
where $A_{1}$, $B_{1}$, and $C_{1}$ are as defined in the statement of the Corollary. Initial conditions and
multiplicities of periodic solutions are unchanged under this transformation. Now the result follows from Theorem B.
\newtheorem*{r3}{Remark}
\begin{r3}
If the differential equation has a linear term, then the equations for $d_{k}$ are more complicated. Instead, we will
have linear differential equations in which the right-hand side depends also on $d_{k}$. We consider such cases in the
next section.
\end{r3}
\section{Proof of Theorems C and D}
\begin{proof}
\emph{(Theorem C):}\\
 We follow the procedure of the last section. In this case, we obtain a sequence of linear
differential equations.
\[\dot{d}_{1}=A d_{1},\]
\[\dot{d}_{k}=A d_{k},\,\, 2 \leq k \leq (N-1),\]
and
\[\dot{d}_{N}=(d_{1})^{N}+A d_{N} \]
together with the initial conditions
\[d_{1}(0)=1;\,d_{k}(0)=0,\, 2 \leq k \leq N. \]
Solving these initial value problems, gives
\[d_{1}(t)=e^{\int_{0}^{t} A(s) ds},\]
\[d_{k}(t)\equiv 0,\,\, 2 \leq k \leq (N-1), \]
and
\[d_{N}(t) = \int_{0}^{t} A(s)ds \int_{0}^{t} e^{(N-1)\int_{0}^{s} A(u)
du} ds.\] The origin is of multiplicity $N$ if and only if $d_{1}(\omega)=1$ and $d_{N}(\omega) \neq 0$. These two
conditions give
\[ e^{\int_{0}^{\omega} A(t) dt }=1, \]
and \[ \int_{0}^{\omega} e^{(N-1)\int_{0}^{t} A(s) ds} dt \neq 0.\] If these two integrals vanish then the solution
$z=0$ is a center. This follows from the general solution of this equation. The general solution is given by
\[ w(t)=e^{(1-N) \int_{0}^{t} A(s) ds} [(1-N) \int_{0}^{t} e^{(N-1)
\int_{0}^{s} A(u) du} ds + C], \] where $w=z^{1-N}$ and $C$ is a constant. Since $w(0)=w(\omega)$, all the solutions
are periodic. Hence, $z=0$ is a center.
\end{proof}

It is clear that the conditions for a center are not satisfied by any real function $A(t)$. We give two equations that
have centers at the origin. In the first equation, we have
\[A(t)=C+ \frac{2}{N-1} \tan (t+c_{0}) .\]
Since the zeros of the complex function $\cos$ are real, the function $A(t)$ is a periodic function of period $2 \pi$,
when $c_{0}$ is a complex non-real number. The origin is of multiplicity one when
\[ e^{\int_{0}^{2 \pi} [C+ \frac{2}{N-1}\tan(t+c_{0}) ]dt} =1.\]
This gives
\[ e^{ 2 \pi C} e^{\frac{1}{N-1}
[\ln(\sec^{2}(2\pi+c_{0}))-\ln(\sec^{2}(c_{0}))]} = 1.\] Since $\sec(t+c_{0}) $ is of period $2 \pi$, the condition
becomes $e^{2 \pi C} =1$. Hence, $C=pi$ for an integer $p$. The second condition is
\[ \int_{0}^{2\pi} e^{(N-1)pit} \sec^{2}(t+c_{0}) dt=0. \]
If $(N-1)p$ is odd, then a change of variables $t \mapsto t-\pi$ gives
\[ \int_{\pi}^{2\pi} e^{(N-1)pit} \sec^{2}(t+c_{0}) dt=-\int_{0}^{\pi} e^{(N-1)pit} \sec^{2}(t+c_{0}) dt=0. \]
Therefore,
\[ \int_{0}^{2\pi} e^{(N-1)pit} \sec^{2}(t+c_{0}) dt=0. \]
 The conditions for a center given in Theorem C are
satisfied. Hence, all the solutions $z(t,c)$, with $c$ is in a neighborhood of $0$, are $2\pi$-periodic. This proves
the first part of Corollary C.

The second part in Corollary C, with $N=3$, is similar to the result of \cite{A6}. In fact, it was shown in \cite{A6}
that $z=0$ is a center for the equation
\[ \dot{z} = z^{3} +(C_{0} + C_{1} e^{-it} + C_{2} e^{it})z
\]
if and only if $C_{0}=ki$ for an integer $k$, $C_{1} C_{2} \neq 0$, and $4i \sqrt{C_{1}C_{2}}$ is a zero for
$J_{2|k|}$.

Now, we prove the second part of Corollary C. The first condition for a center in Theorem C is
\[ e^{\int_{0}^{2\pi} (C - xi \sin t)dt}=e^{2C\pi}.\]
But, $e^{2C \pi}=1$ if and only if $C=pi$, for an integer $p$.
 The
second condition for a center becomes
\[\int_{0}^{2\pi}e^{pit-xi \sin t} dt = \int_{0}^{2\pi} \cos(pt - x
\sin t) dt + i \int_{0}^{2\pi} \sin(pt - x \sin t) dt.\] The imaginary part of this quantity is zero; it is an integral
of an odd $2\pi-$periodic function over the interval $[0,2\pi]$. The real part is $2\int_{0}^{\pi} \cos(pt - x \sin t)
dt$. This integral is the integral form of the Bessel's function. Therefore The real part is zero if $x$ is a zero of
the Bessel's function $J_{p}$. We recall that the integral form and power series expansion of $J_{p}$ are defined by
\[J_{p}(x)= \int_{0}^{\pi} \cos(pt - x \sin t) dt= \sum_{j=0}^{\infty}
(-1)^{j}\frac{(\frac{x}{2})^{p+2j}}{j!(p+j)!}.\]

The power series representation of Bessel's function is used in \cite{A6}. Consequently, their proof is much longer
than our proof.

\begin{proof}
\emph{(Theorem D:)}\\
 With the change of variables  $z \mapsto e^{\int_{0}^{t}q(s) ds} z$, the equation becomes
\[ \dot{z}=\sum_{k=2}^{N} p'(t)\,f_{k}(p(t))\,z^{k}. \]
From the expansion of Section 3, we can show inductively that the coefficients $d_{k}$ are functions of $p(t)$. This
implies that $z(t,c)$ is a function of $p(t)$. Therefore, in a neighborhood of the origin, all solutions are
$\omega-$periodic and $z=0$ is a center.
\end{proof}
\newtheorem*{r7}{Remark}
\begin{r7}
If the coefficient of $z^{N}$ in the statement of Theorem D satisfies the condition
\[ p'(t)\,f_{N}(p(t))=e^{(N-1)\,\int_{0}^{t}q(s)ds} \]
then the coefficient of $z^{N}$ equals $1$. This condition is satisfied when $p(t)=\frac{e^{2 \pi i t}}{2 \pi i}$,
$f_{N} \equiv 1$, and $q=2 \pi i$. The period in this case is $N-1$.
\end{r7}
\section{The Pugh problem}
The use of computer algebra has led to significant progress in the investigation of the properties of polynomial
differential systems. In this section, we describe an application of computer algebra to find the maximum possible
multiplicity of periodic solutions of polynomial differential equations. In particular, we present solutions to the
local Pugh problem \cite{A13} and Shahshahani conjecture \cite{A12}. We mention that Pugh problem was listed as a part
of problem 13 in Steve Smale list of 18 open problems for the next century. This problem was considered as a version of
Hilbert sixteenth problem. Hilbert sixteenth problem is to estimate the number of limit cycles of polynomial
two-dimensional systems. Research related to Hilbert sixteenth problem has derived enormous benefit from the
availability of computer algebra as is demonstrated by, for instance, \cite{A11}, and \cite{A14}.

Consider the differential equation
\begin{equation}
\dot{z}=x^{4}+A(t)x^{3}+B(t)x^{2}
\end{equation}

We show that this local problem can be solved using the method of Groebner bases. We give an automatic means of finding
the maximum possible multiplicity of a periodic solution. The algorithm involves computing Groebner bases. Computer
algebra systems, such as Maple, can be used to implement this algorithm. If the equation does not have a complex
center, then the local Pugh problem is solvable by our procedure. The algorithm for computing the maximum possible
multiplicity is then described in this section. We apply the algorithm for equations in which the coefficients $A(t)$
and $B(t)$ are polynomial functions in $t$, and in $\cos t$ and $\sin t$. The case $N=4$ is considered. However, the
method works for any $N$. In the case that the coefficients are polynomial functions in $t$, we show that the maximum
possible multiplicity is the same whether the coefficients are complex or are restricted to be real. When the
coefficients are trigonometric polynomials, cases are described where the maximums are not equal. Moreover, we show
that the origin is not a complex center in each of the equations considered; these results provide evidences to support
our conjecture.

We follow the same procedure of Section 3. For equation (5.1), $d_{1}(t)\equiv1$ and the equations satisfied by the
$d_{n}(t)$ (for $n>1$) are
\begin{equation}
\dot{d}_{n}=\sum_{i+j+k+l=n}^{}d_{i}d_{j}d_{k}d_{l}+A\, \sum_{i+j+k=n}^{}d_{i}d_{j}d_{k}+ B\,
\sum_{i+j=n}^{}d_{i}d_{j}.
\end{equation}
These equations were integrated by parts repeatedly and the formulae for $d_{n}$, with $n\leq 8$, are given in
\cite{A3}. The calculations become extremely complicated as $n$ increases. It is impossible to accomplish these
computations by hand except in the simplest cases.

The next step, in computing the multiplicity, is to consider the quantities
\begin{equation*}
 \eta_{n}=d_{n}(\omega)
 \end{equation*}
  Note that$\eta_{n}$ is a polynomial function in the coefficients of the
polynomials $A$ and $B$. The multiplicity of the origin is $k$ if
\[\eta_{2}=\eta_{3}=\cdots =\eta_{k-1}=0,\eta_{k} \neq 0.\]
We write $\mu_{max}(\mathcal{C})$ for the maximum possible multiplicity of $z=0$ for equations in a class
$\mathcal{C}$. For the class of equations in which the coefficients are polynomial functions of degree $m$, Pugh
problem is to find $\mu_{max}$ in terms of $m$.

Now, we give the formulae of $\eta_{n}$ with $n\leq 4$, for equation {5.1}. We use a tilde over a function to denote
its indefinite integral:
\[
\widetilde{f}(t)=\int_{0}^{t}f(s)ds
\]
\newtheorem*{T5}{Proposition}
\begin{T5}
For equation (5.1), the quantities $\eta_{2},\eta_{3}$, and $\eta_{4}$ are as follows.
\[\eta_{2}=\int_{0}^{\omega}B dt \]
\[\eta_{3}=\int_{0}^{\omega}A dt\]
\[\eta_{4}=\omega+ \int_{0}^{\omega}(\widetilde{B}A)dt\]
\end{T5}
\begin{proof}
The formulae were obtained by solving equations (5.2) recursively. The computations leading to the formulae proceed by
sequences of judiciously chosen integrations by parts; for any functions $f$ and $g$, we make use of the identity
\[\widetilde{f\widetilde{g}}=\widetilde{f}\widetilde{g}-\widetilde{\widetilde{f}g}.\]
These are elementary though together they form a complicated web. To obtain $\eta_{n}$, we reduce $d_{n}(\omega)$
modulo the ideal generated by $d_{2}(\omega),\cdots,d_{n-1}(\omega)$.
\end{proof}

\noindent
 \textbf{The Algorithm}\\
 We shall write $\mu_{max}(real)$ when the coefficients of the polynomials $A$ and $B$ are real numbers, and
$\mu_{max}(complex)$ when the coefficients are complex numbers. It follows from the result in \cite{A1}, that
$\mu_{max}(real) < \infty$. We call the set of equations that have this maximum multiplicity, the maximum variety,
$V_{max}$. Similarly, we define $V_{max}(real)$ and $V_{max}(complex)$.

To find $\mu_{max}$, first we integrate recursively to compute the functions $d_{n}(t)$. Then we consider the
expressions $d_{n}(\omega)$, which are polynomial functions in the coefficients of $A$ and $B$. To obtain $\eta{n}$ we
reduce $d_{n}(\omega)$ modulo the ideal generated by $d_{2}(\omega),\cdots,d_{n-1}$. We stop until the system
$\eta_{2}=\eta_{3}=\cdots=\eta_{k}=0$ has no real solutions. In this case $\mu_{max}=k$. From the theory of Groebner
bases, the Groebner basis of the ideal $\langle \eta_{2},\eta_{3},\cdots,\eta_{k} \rangle$ is $\langle 1 \rangle$ if
and only if the system $\eta_{2}=\eta_{3}=\cdots=\eta_{k}=0$ has no complex solutions. So, we have to verify that this
maximum multiplicity can be attained by certain real values of coefficients. The procedure gives an upper bound for
$\mu_{max}$. However, for the equations in which the coefficients are polynomial functions of $t$ that we will
consider, $\mu_{max}(real) = \mu_{max}(complex) < \infty$.

The algorithm for computing $\mu_{max}$ can be summarized as follows:

\noindent
 $\bullet$ Input: functions
$A(t)$ and $B(t)$ which are polynomials in $t$, or in $\cos t$ and $\sin t$. \\ $\bullet$ Integrate to compute
$d_{n}(t)$.\\ $\bullet$ Compute $d_{n}(\omega)$.\\ $\bullet$ Find $\eta_{n}$ by reducing $ d_{n}(\omega)$ modulo
$\langle \eta_{2},\eta_{3},\cdots,\eta_{n-1} \rangle$.\\ $\bullet$ Stop when the Groebner basis $\langle
\eta_{2},\eta_{3},\cdots,\eta_{k}\rangle$ is $\langle 1 \rangle$.\\ $\bullet$ Output: $\mu_{max}=k$. For the details
related to Groebner bases, we refer to \cite{A7}. We use Maple8 to compute Groebner bases.

Now, we consider the first class of coefficients. Here, we assume that $\omega=1$.

Let $B(t) = C_{1}+at+ bt^{2}$ and $ A(t)=C_{2}+ct+dt^{2}+et^{3}.$ If $\eta_{2}=\eta_{3}=0$, then
\[ 12C_{1}+6a+4b+3c=0, 12C_{2}+6d+4e+3f=0.\]
This class of equations was used in \cite{A2} to construct equations with $10$ real periodic solutions. We substitute
the values of $C_{1}$ and $C_{2}$ from these equations and then compute the Groebner basis of the ideal $\langle
\eta_{4},\eta_{5},\cdots,\eta_{10} \rangle $. This basis is $\langle 1 \rangle $. Therefore, the maximum possible
multiplicity is $10$. To find the equations that have this maximum multiplicity, we compute the Groebner basis $\langle
\eta_{4},\eta_{5},\cdots,\eta_{9}\rangle$. This set of equations has $15$ solutions, counting multiplicity, and at
least one of the solutions is real. We summarize the result for this class in the following lemma.
\newtheorem*{T6}{Lemma 1}
\begin{T6}
If $A(t)$ and $B(t)$ are polynomial functions in $t$ of degrees $3$ and $2$, respectively, then
$\mu_{max}(real)=\mu_{max}(complex)=10$. Moreover, $V_{max}$ is a zero-dimensional ideal.
\end{T6}
The Groebner basis is computed with the graded reverse lexicographic term order and with respect to the list
$[b,c,d,a,e]$; this term order usually gives more compact Groebner basis. Then, the basis is changed to the
lexicographic term order, which is the most suitable to eliminate variables from a set of equations. The Groebner basis
is given at the end of this section. The last equation in this basis is a polynomial in $e$ of degree $15$; it has at
least one real solution. The other variables are given explicitly as functions of $e$.

The next class has the coefficients
\[ B(t)= C_{1}+at,A(t)=C_{2}+bt+ct^{2}+dt^{3}+et^{4}+ft^{5}. \]
 Let $2C_{1}+a=0$ and $60C_{2}+30b+20c+15d+12e+10f=0$ for $\eta_{2}=\eta_{3}=0$. We compute the Groebner basis
 as in the first case. This gives $\langle \eta_{4},\eta_{5},\cdots,\eta_{10}\rangle=\langle 1 \rangle$.
 Moreover,
 $\langle \eta_{4},\eta_{5},\cdots,\eta_{9}\rangle$ has a polynomial equation in $f$ of degree $6$, which has
 two real solutions and four complex non-real solutions. The other variables are given as functions of $f$.
 The result for this class are given in the following lemma.
\newtheorem*{T7}{Lemma 2}
\begin{T7}
If $A(t)$ and $B(t)$ are polynomial functions in $t$ of degrees $5$ and $1$, respectively, then
$\mu_{max}(real)=\mu_{max}(complex)=10$. Moreover, $V_{max}$ is a zero-dimensional ideal.
\end{T7}
The Groebner basis is given at the end of the section.

In the last class of polynomials in $t$, we let
\[ A(t) = C_{1}+bt+ ct^{2} + d t^{3} , B(t)=C_{2}+at.\]
If $\eta_{2}=\eta_{3}=0$, then
\[ 12C_{1}+6b+cb+3d=0, 2C_{2}+a=0.\]
We substitute the values of $C_{1}$ and $C_{2}$ from these equations and then compute the Groebner basis of the ideal
$\langle \eta_{4},\eta_{5},\cdots,\eta_{8}\rangle$. This basis is $\langle 1 \rangle$. Therefore, the maximum possible
multiplicity is $8$. To find the equations that have this maximum multiplicity, we compute the Groebner basis $\langle
\eta_{4},\eta_{5},\cdots,\eta_{7}\rangle$. The last equation, which is a polynomial in $a$, has two real solutions and
four complex non-real solutions. The other variables are given in terms of $a$. Again, the basis is given after the
statement of Lemma 3.
\newtheorem*{T8}{Lemma 3}
\begin{T8}
If $A(t)$ and $B(t)$ are polynomial functions in $t$ of degrees $3$ and $1$, respectively, then
$\mu_{max}(real)=\mu_{max}(complex)=8$. Moreover, $V_{max}$ is a zero-dimensional ideal.
\end{T8}
\begin{gather*}
\langle 8695641600 \,b-6086949120\, d-773773\, a^{5}+3151791360\, a^{2},\\ 8695641600 c+13043462400
d+773773\,a^{5}-3151791360\, a^{2},\\ -50295245\, a^{4}+37439568 \,d^{2}+118513886400\, a,\\ -120772800 \,d+70343 \,d
a^{3}, 773773\, a^{6}-3151791360\, a^{3}+3130430976000 \rangle
\end{gather*}

Now, we consider classes of coefficients which are polynomial functions in $\cos t$ and $\sin t$. Here, we take
$\omega=2 \pi$.
\newtheorem*{T9}{Lemma 4}
\begin{T9}
If $A(t)$ and $B(t)$ are homogeneous polynomial functions in $\cos t$ and $\sin t$ of degree $1$ or $2$, then
$\mu_{max}(real)=6$ and $\mu_{max}(complex)=7$. The sets $V_{max}(real)$ and $V_{max}(complex)$ are not
zero-dimensional ideals.
\end{T9}
In the case $A(t)=c \cos t + d \sin t$, and $B(t)=a \cos t + b \sin t$. It is clear that $\eta_{2}=\eta_{3}=0$. From
Maple, we have:
\[ \langle \eta_{4},\eta_{5},\eta_{6}\rangle= \langle a^{2}+b^{2},bc-ad-2  \rangle. \]
This set has complex solutions but does not have a real solution. Moreover
\[ \langle \eta_{4},\eta_{5},\eta_{6},\eta_{7}\rangle = \langle 1 \rangle .\]
For the other case $A(t)=c \cos^{2} t + d \cos t \sin t + c_{1} \sin^{2} t$ and $B(t)=a \cos^{2} t + b \cos t \sin t +
a_{1} \sin^{2} t$. We have $\eta_{2}=a+a_{1}$ and $\eta_{3}=c+c_{1}$. If $a_{1}=-a$ and $c_{1}=-c$, then Maple gives
\[ \langle \eta_{4},\eta_{5},\eta_{6}\rangle= \langle 4a^{2}+b^{2},bc-ad-8  \rangle, \]
and
\[ \langle \eta_{4},\eta_{5},\eta_{6},\eta_{7}\rangle = \langle 1  \rangle. \]
A similar argument proves the Lemma.

If one of the coefficients $A(t)$ and $B(t)$ contains only terms with of even degrees and the other coefficient
contains only terms of odd degrees, then it follows from the Proposition that $\mu_{max}(real)=\mu_{max}(complex)=4$.
When the coefficients are homogeneous polynomials of degree $3$, the Groebner basis becomes very large.
\newtheorem*{T10}{Lemma 5}
\begin{T10}
Let $A(t)$ be a homogeneous polynomial of degree $3$, and $B(t)$ is a homogeneous polynomial of degree $1$. The
solution $z=0$ is not a center. Moreover, $\mu_{max}(real)=8$, and $\mu_{max}(complex)=10$.
\end{T10}
We take $B(t)= a \cos t + b \sin t$. Using the identities $\sin^{3} t = \sin t(1- \cos^{2} t)$ and
 $\cos^{3}(t)=\cos t(1-\sin^{2} t)$, we can write a homogeneous polynomial of degree $3$ in the following form:
\[ A(t)=c \cos t + d \sin t + e \cos t \sin^{2} t + f \sin t \cos^{2} t.\]
The computer output is given at the end of this section. The Groebner basis of $\langle
\eta_{4},\eta_{5},\eta_{6},\eta_{7} \rangle$ has a real solution; we can take $a=d=1,c=f=0, e^{2}=432$, and $36\,b=e$.
On the other hand, $\langle \eta_{4},\cdots,\eta_{8}\rangle$ does not have a real solution but $\langle
\eta_{4},\cdots,\eta_{9} \rangle \neq \langle 1 \rangle$, and $\langle \eta_{4},\cdots, \eta_{10} \rangle =\langle 1
\rangle $. On the other hand, $\eta_{8}$ has the form
\[ \eta_{8}=\pi\,((5/8)\,a^{4}+(5/8)\,b^{4}+(5/4)\,a^{2}\,b^{2} +(1/96)\,e^{2} +(1/96)\,f^{2}) .\]
 It is clear that $\eta_{8} \neq 0$,
when the coefficients are restricted to be real numbers.
\newpage
\emph{Maple output for Lemma 1.}\\
 $\langle \eta_{4},\cdots,\eta_{9} \rangle=$
\begin{gather*}
\langle 5077631499817345768910374214846690690808088946119335814119700628085731\\ 205306414251948895671156736000
\,a+\\76568765938297226294215511996032531264938846500344266059142283505\,e^{14}-\\
4136638384639181060935454884796563027263107079748408223299514765436040 5024 \,e^{11}\\
+9607203693297398146546053698859951537568361672283979690465766480307919 \\005421113600 \,e^{8}-\\
1077298213530002385272086712469698588685019454221289878562607166668317\\ 844870406197314322432 \,e^{5}\\
+5120940699074104547488742022229220044464651558296140078787235311492364 \\0199838438141019890846269440 \,e^{2},\\
5077631499817345768910374214846690690808088946119335814119700628085731\\ 205306414251948895671156736000
\,b-\\76568765938297226294215511996032531264938846500344266059142283505\,e^{14}+\\
4136638384639181060935454884796563027263107079748408223299514765436040 5024 \,e^{11}\\
-9607203693297398146546053698859951537568361672283979690465766480307919\\ 005421113600 \,e^{8}+\\
1077298213530002385272086712469698588685019454221289878562607166668317\\ 844870406197314322432 \,e^{5}-\\
5120940699074104547488742022229220044464651558296140078787235311492364\\ 0199838438141019890846269440 \,e^{2},\\
2542195853515507102698935612213344207862818877038578400901689200003189 \\8438696161676309299200
\,c\\-8114451463043258106360200365311543313307196783429471460595\,e^{13}+\\
4874645150831477756390000057954655294888051734933883104748830146656\, e^{10}-\\
1175665268311675043075984579090225104206445216121635282751902996942202\\ 131200 \,e^{7}+\\
1325591495257443095378482025935592779323248112155531380849978417494633 \\00399322103808 \,e^{4}-\\
2197311395713698234662544085226573648290742502622945660686244687349190\\ 8621199707934252072960
\,e,\\2\,d+3\,e,\\21063213975871712330828163101314609697689436160\,e^{3}\\
+3881552635098624688463371219200\,e^{9}-16555584308626955352096\,e^{12}\\
-438195979208967285885778234437995593728\,e^{6}+29908464432145
\,e^{15}\\-3407873245453536523174061109544883049290583244800000 \rangle
\end{gather*}
\emph{Maple output for Lemma 2.}\\ $\langle \eta_{4},\cdots,\eta_{9} \rangle =$
\begin{gather*}
 \langle 12699270921091068840199720844433641567807\\
 091515760\,b+\\
2539854184218213768039944168886728313561418303152\,e+\\
 423247429138101258858306608404799875\,f^{4}+\\
4341516136354618897336295790080428452939652335036\,f,\\
 496186409615279857004140496974208217545576225280\,c-\\
595423691538335828404968596369049861054691470336\,e-\\
 38477039012554659896209691673163625\,f^{4}-\\
1028512938658225137299402571203191057331286500736\,f,\\
 38477039012554659896209691673163625\,f^{4}+\\
744279614422919785506210745461312326318364337920\,d+\\
 1488559228845839571012421490922624652636728675840\,e+\\
2020885757888784851307683565151607492422438951296\,f,-\\
 76954078025109319792419383346327250\,ef^{4}-\\
192385195062773299481048458365818125\,f^{5}+\\
 18755846283457578594756510785625070623222781315584000\,a-\\
6849063961008228076202649414049259748983323392\,e f-\\
 17122659902520570190506623535123149372458308480\,f^{2},\\
4554469328710172509947246108104831527932968138654523985\\ 40570112\,e^{2}\\
+2277234664355086254973623054052415763966484069327261992\\ 702850560\,e f\\
28452809860130268863272755268929911307472870565806959073\\ 22835840\,f^{2}\\
-49539000048633777493317082388741947622188883125\,f^{5},\\
 38477039012554659896209691673163625\,f^{6}+\\
3984808136948929447103185415887785672672861696\,f^{3}-\\ 14218838340357434859762303459719245720157600157874298880000
\rangle
\end{gather*}
\emph{Maple output for Lemma 5.}\\ $ \langle a-1,d-1,c,f,\eta_{4},\cdots,\eta_{7} \rangle =
[-432+{e}^{2},a-1,36\,b-e,c,d-1,f]$\\
 $ \langle
\eta_{4},\cdots, \eta_{7} \rangle=$
\begin{gather*}
[f{a}^{4}+2\,{b}^{2}f{a}^{2}+{b}^{4}f-24\,{a}^{3}+72\,{b}^{2}a,8\,cf{a
}^{3}+16\,{a}^{2}fbd+8\,{b}^{3}df+7\,{a}^{2}{f}^{2}b+3\,{b}^{3}{f}^{2}\\
-192\,{a}^{2}c+576\,abd-16\,{b}^{2}e+168\,afb-384\,b, 24\,f{a}^{2}{c}^{
2}+48\,{d}^{2}f{a}^{2}+24\,{d}^{2}f{b}^{2}+42\,{a}^{2}d{f}^{2}\\ +18\,d{f }^{2}{b}^{2}+
9\,{f}^{3}{a}^{2}+3\,{f}^{3}{b}^{2}-576\,a{c}^{2}+1728\,{ d}^{2}a-96\,bde\\
+4\,{e}^{2}a+1008\,adf-18\,fbe+142\,{f}^{2}a-2304\,d- 912\,f, 6912\,fa{c}^{3}+13824\,{d}^{2}fac+6912\,{d}^{3}fb\\
-1536\,d{e}^{2 }bf +4\,{e}^{3}fa+12960\,dc{f}^{2}a+9504\,{d}^{2}b{f}^{2}-406\,{f}^{2}b {e}^{2}+ 2916\,ac{f}^{3}\\
+2676\,d{f}^{3}b+31\,e{f}^{3}a+161\,{f}^{4}b
 -165888\,{c}^{3}+497664\,{d}^{2}c-82944\,e{c}^{2}\\
 +165888\,e{d}^{2}-
10368\,{e}^{2}c -384\,{e}^{3}+490752\,dcf+130176\,def+105120\,c{f}^{2} +23664\,e{f}^{2},f{a}^{3}\\
+3\,a{b}^{2}f-24\,{a}^{2}+48\,{b}^{2}+2\,{b}^{ 3}e, -297\,ec{f}^{2}-297\,d{e}^{2}f-1296\,df{c}^{2}-1296\,e{d}^{2}c\\
-432 \,{f}^{2}{c}^{2}-1296\,defc+432\,e{c}^{3}+216\,{e}^{2}{c}^{2}+27\,{e}^ {3}c
-432\,{e}^{2}{d}^{2}+27\,d{f}^{3}-52\,{e}^{2}{f}^{2}\\ +{e}^{4}+{f}^{
4}+432\,{d}^{3}f+216\,{d}^{2}{f}^{2},8\,{b}^{2}de+4\,{a}^{2}fc+12\,abd f +3\,{b}^{2}ef+5\,a{f}^{2}b-96\,ac+192\,bd\\
-8\,ae+64\,bf,96\,e{d}^{2}b- 2\,{e}^{3}b+48\,fa{c}^{2}+144\,{d}^{2}fa+72\,debf-{e}^{2}fa+120\,ad{f} ^{2}\\
+10\,{f}^{2}be+23\,{f}^{3}a-1152\,{c}^{2}+2304\,{d}^{2}-192\,ce+16
\,{e}^{2}+1536\,df+280\,{f}^{2},{b}^{2}{e}^{2}+{f}^{2}{b}^{2}+12\,be\\ +12\,fa
-288,-2\,afb+{a}^{2}e-{b}^{2}e-24\,b,-f{b}^{2}+2\,aeb+f{a}^{2}- 24\,a,{e}^{2}a+12\,eac-12\,bde\\ -24\,adf-3\,fbe
-8\,{f}^{2}a-288\,d-96\,f , 12\,fac+2\,{e}^{2}b-288\,c+24\,dea+9\,efa-48\,e-12\,dfb-b{f}^{2},\\ -8+b e-fa-4\,ad+4\,cb]
\end{gather*}
$ \langle \eta_{4},\cdots,\eta_{9} \rangle =$
\begin{gather*}
[2\,{e}^{3}c-2\,d{e}^{2}f+2\,ec{f}^{2}-{e}^{2}{f}^{2}-2\,d{f}^{3}-{f}^
{4}-11520\,{a}^{2}-23040\,{b}^{2},2\,{e}^{3}d+2\,f{e}^{2}c\\+f{e}^{3}+2 \,{f}^{2}ed\\
+2\,{f}^{3}c+{f}^{3}e-11520\,ab,{e}^{4}+2\,{e}^{2}{f}^{2}+{ f}^{4}+34560\,{a}^{2}+34560\,{b}^{2},1440\,{a}^{3}\\
-9\,{e}^{2}d-2\,f{e} ^{2} -9\,d{f}^{2}-2\,{f}^{3},1440\,{a}^{2}b+3\,{e}^{2}c+{e}^{3}+3\,c{f}
^{2}+e{f}^{2},1440\,{b}^{2}a-3\,{e}^{2}d\\ -f{e}^{2}-3\,d{f}^{2}-{f}^{3},
1440\,{b}^{3}+9\,{e}^{2}c+2\,{e}^{3}+9\,c{f}^{2}+2\,e{f}^{2},-2\,afb+{ a}^{2}e-{b}^{2}e-24\,b,\\
-f{b}^{2}+2\,aeb+f{a}^{2}-24\,a, {e}^{2}a+{f}^{2 }a+288\,d+72\,f, {e}^{2}b+b{f}^{2}-288\,c-72\,e,\\
24\,ac+5\,ae+bf,24\,cb +7\,be +fa-24,72\,{c}^{2}+39\,ce+5\,{e}^{2}-3\,df-{f}^{2},\\
 be+24\,ad+7\,f a+24,
24\,bd+ae+5\,bf,24\,cd+5\,ed+5\,cf+ef, 72\,{d}^{2}-3\,ce+39\,df+5 \,{f}^{2}-{e}^{2}]
\end{gather*}
\bibliographystyle{amsplain}

\begin{thebibliography}{999}
\bibitem{A1} M.A.M. Alwash; \emph{Periodic solutions of Abel differential
equation}, J. Math. Anal. Appl., 329(2007)1161-1169.
\bibitem{A2} M.A.M. Alwash; \emph{Periodic solutions of a quartic
differential equations and Groebner bases}, J. Comp. Appl. Math., 75(1996)67-76.
\bibitem{A3} M.A.M. Alwash and N.G. Lloyd; \emph{Periodic solutions
of a quartic nonautonomous equation}, Nonlinear Analysis, 11(1987)809-820.
\bibitem{A4} M.A.M. Alwash and N.G. Lloyd; \emph{Non-autonomous
equations related to polynomial two-dimensional systems}, Proc. Royal Soc. Edinburgh, 105(1987)129-152.
\bibitem{A5} M. Briskin and Y. Yomdin; \emph{Tangential version of
Hilbert 16th problem for the Abel equation}, Mosc. Math. J., 5(2005)23-53.
\bibitem{A6} A. Cima, A. Gasull, and F. Manosas; \emph{ Periodic
orbits in complex Abel equations}, J. Diff. Eqns., 232(2007) 314-328.
\bibitem{A7} D.A. Cox and B. Sturmfels; \emph{Applications of Computational Algebraic Geometry}, American
Mathematical Society(1998).
\bibitem{A8} Y. Ilyashenko; \emph{Hilbert-type numbers for genralized Abel
equations, growth and zeros of holomorphic functions}, Nonlinearity, 13(2000)1337-1342.
\bibitem{A9} A. Lins Neto; \emph{On the number of solutions of the equation
$\frac{dx}{dt}=\sum_{j=0}^{n}a_{j}(t)x^{j}, 0 \leq t \leq 1,$ for which $x(0)=x(1)$}, Invent. Math., 59(1980)67-76.
\bibitem{A10} N.G. Lloyd; \emph{The number of periodic solutions of the
equation $\dot{z}=z^{N}+P_{1}(t)z^{N-1}+\cdots+P_{N}(t)$}, Proc. London Math. Soc., 27(1973)667-700.
\bibitem{A11} N.G. Lloyd and J.M. Pearson; \emph{Symmetry in planar dynamical systems}, J. Symbolic Computation,
33(2002)357-366.
\bibitem{A12} S. Shahshahani; \emph{Periodic solutions of polynomial
first order differential equations}, Nonlinear Analysis, 5(1981)157-165.
\bibitem{A13} S. Smale; \emph{Mathematical problems for the next
century}, Mathematics: Frontiers and Perspectives, AMS(2000), 271-294.
\bibitem{A14} D. Wang; \emph{Polynomial systems from certain differential equations}, J. Symbolic Computation,
28(1999)303-315.
\bibitem{A15} Y. Yomdin; \emph{The center problem for the Abel equations, compositions of functions, and moment
conditions},
Mosc. Math. J., 3(2003)1167-1195.
\end{thebibliography}

\end{document}